\documentclass[a4paper,12pt]{article}
\usepackage[top=2.5cm,bottom=2.5cm,left=2.5cm,right=2.5cm]{geometry}
\usepackage{graphicx}
\usepackage{amsmath,amsthm,amssymb,lineno}
\usepackage{latexsym}
\usepackage{graphicx,color, booktabs,multirow}
\usepackage{tikz}
\usepackage{graphicx,booktabs,multirow,mathrsfs}
\usepackage{appendix}
\usepackage{yhmath}
\usepackage{setspace}
\usepackage{soul}
\usepackage{enumitem}
\newtheorem{theo}{Theorem} %[section] %
\newtheorem{lemma}[theo]{Lemma}
\newtheorem{conj}{Conjecture}

\newtheorem{prop}[theo]{Proposition}
\newtheorem{prob}{Problem}
\newtheorem{claim}{Claim}
\usetikzlibrary{backgrounds}
\usetikzlibrary{arrows}
\usetikzlibrary{shapes,shapes.geometric,shapes.misc}
\pgfdeclarelayer{edgelayer}
\pgfdeclarelayer{nodelayer}
\pgfsetlayers{background,edgelayer,nodelayer,main}
\tikzstyle{none}=[inner sep=0mm]
\tikzstyle{every loop}=[]

\tikzstyle{new style 0}=[fill=black, draw=black, shape=circle]
\tikzstyle{red style 1}=[fill=red, draw=black, shape=circle]
\tikzstyle{blue style 2}=[fill=blue, draw=black, shape=circle]
\tikzstyle{white style 4}=[fill=white, draw=black, shape=circle]
\tikzstyle{bklack style 5}=[fill=black, draw=black, shape=rectangle]
\tikzstyle{red style 3}=[fill=red, draw=black, shape=rectangle]
\tikzstyle{yellow style 7}=[fill=yellow, draw=black, shape=rectangle]
\tikzstyle{new style 8}=[fill={rgb,255: red,0; green,132; blue,0}, draw={rgb,255: red,0; green,131; blue,0}, shape=circle]

\tikzstyle{new edge style 0}=[-]
\tikzstyle{new edge style 1}=[-, draw=red]
\tikzstyle{new edge style 2}=[-, draw=blue]
\tikzstyle{new edge style 3}=[-, draw={rgb,255: red,0; green,156; blue,0}]

\newcounter{countcase}

\def \proof {\noindent {\bf Proof}. }
\def \proofend {\hfill $\Box$  \setcounter{countcase}{0}}

\newcommand {\red} {\textcolor{red}}
\newcommand {\green} {\textcolor{green}}
\newcommand {\blue} {\textcolor{blue}}
\newcommand {\grey} {\textcolor{grey}}
\def \F {{\mathscr F}}

\newcommand\cro[1]
{
{#1}^{\times}
}

\newcommand\equ[2]
{
\begin{equation}\label{#1}
#2
\end{equation}
}

\def \usecolour
{\newcommand {\rered}{\red}
\newcommand {\reblue} {\blue}
\newcommand {\regreen} {\green}
\newcommand {\regrey} {\grey}
}

\def \nocolour
{\renewcommand {\rered}{}
\renewcommand {\reblue}{}
\renewcommand {\regreen}{}
\renewcommand {\grey}{}
}
\usecolour

\nocolour    %for the final
\newcommand \rebibitem[1]
{
\bibitem{#1} %\red{\bf ***#1***}
}

\newcommand \notwant[1]
{\newpage
\red{\bf *****SOME IDEAS*****}

\vspace{0.6 cm}

#1\newpage
}

\renewcommand \notwant[1]{}

\begin{document}
%\linenumbers

\title{\bf On the size of matchings in 1-planar graph with high minimum degree\thanks{This work is supported by MOE-LCSM, School of Mathematics and Statistics,
Hunan Normal University, Changsha, Hunan, %410081,
China  and Hunan Provincial Natural Science Foundation of
China (No. 2021JJ30169).}}
\author{ Yuanqiu Huang \\
{\footnotesize  Department of Mathematics, Hunan Normal University, Changsha 410081, P.R.China} \\
{\footnotesize hyqq@hunnu.edu.cn}\\
Zhangdong Ouyang\thanks{Corresponding author.}\\
{\footnotesize Department of Mathematics,
Hunan First Normal University,
Changsha 410205, P.R.China} \\
{\footnotesize oymath@163.com }\\
Fengming Dong\\
{\footnotesize National Institute of Education,
Nanyang Technological University, Singapore} \\
{\footnotesize fengming.dong@nie.edu.sg
and donggraph@163.com}\\
}

\date{}
\maketitle

\begin{abstract}
A matching of a graph  is a set of edges  without common end vertex.   A graph  is called 1-planar  if it admits a drawing in the plane such that each edge is crossed at most once.
Recently, Biedl and  Wittnebel proved that every 1-planar graph with minimum  degree 3 and $n\geq 7$ vertices  has a matching  of size at least $\frac{n+12}{7}$, which is tight for some graphs.
They also provided tight lower  bounds  for the sizes of matchings in
1-planar graphs with minimum degree 4 or 5.
In this paper, we show that any 1-planar graph with minimum degree 6 and $n \geq 36$ vertices has a matching of size at least $\frac{3n+4}{7}$, and this lower bound is tight. Our result  confirms a  conjecture posed by Biedl and  Wittnebel.

 \vskip 0.2cm
\noindent
%\textbf{AMS classification}: 05C07, 05C15, 05C50\\
\end{abstract}

\maketitle

%%%%%% THIS PART MUST BE PLACED IMMEDIATELY AFTER THE \maketitle COMMAND
%%%%%% BACK TO ORIGINAL FOOTNOTES
\makeatletter
\renewcommand\@makefnmark%
{\mbox{\textsuperscript{\normalfont\@thefnmark)}}}
\makeatother
%%%%%%

\section{Introduction}

 A {\em drawing} of a graph
$G=(V,E)$ is a mapping $D$ that assigns to each vertex in $V$ a
distinct point in the plane and to each edge $uv$ in $E$ a
continuous arc connecting $D(u)$ and $D(v)$. We often make no
distinction between a graph-theoretical object (such as a vertex, or
an edge) and its drawing. All drawings considered here are
such ones  that no edge crosses itself, no
two edges cross more than once, and no two edges incident with the
same vertex cross. A graph is {\it planar} if it can be drawn in the plane without edge crossings.  A drawing of a graph is 1-{\it planar} if each of its edges is
crossed at most once. If a graph has a 1-planar drawing, then it is
1-{\it planar}. The notion of 1-planarity was introduced in 1965 by
Ringel \cite{GR}, and since then many properties of 1-planar graphs have
been studied (e.g. see the survey paper \cite{SK}). For example, it is known that \cite{HR,IT, JP} any 1-planar graph with $n$ $(\geq 3)$ vertices has at most $4n-8$ edges, and thus has the minimum degree $\leq 7$.

A {\it matching} of a graph is a set of edges without common end vertex.
The study of matchings is one of the oldest and best-studied problems in graph theory, for example, see \cite{LP}. An earlier result,
due to Nishizeki and Baybars \cite{NB},
shows that every simple planar graph with $n\geq X$ vertices has a matching of size at least $Yn+Z$, where $X, Y, Z$ depend on the minimum degree and the connectivity of the graph. In recent years, researchers have investigated the sizes of the matchings graphs that are
``almost" planar.
An interesting example is 1-planar graphs, which  is a generalization  of  planar graphs,
in some sense.
Many papers studying the sizes of matchings in
1-planar graphs have been published,
for example, see
(\cite{B, TB, BK, BW}).
Recently, Biedl and Wittnebel~\cite{BW}
obtained  the following results on
the lower bounds of sizes of matchings in 1-planar graphs
in terms of their minimum degrees.

\begin{theo}[\cite{BW}]\label{th1}
Any $n$-vertex simple 1-planar graph with minimum degree $\delta$ has a matching $M$ of the following size:

   1. $|M|\geq \frac{n+12}{7}$ if $\delta =3$ and $n\geq 7$;

   2. $|M|\geq \frac{n+4}{3}$ if $\delta =4$
and $n\geq 20$; and

   3. $|M|\geq \frac{2n+3}{5}$ if $\delta =5$ and $n\geq 21$.
\end{theo}

The authors in \cite{BW} constructed
$1$-planar graphs which contain matchings
of sizes equal to the lower  bounds in Theorem~\ref{th1}.
As for the minimum degree $\delta=6$,
they also  constructed  1-planar graphs which
\rered{contain} matchings with a maximum size.

\begin{theo}[\cite{BW}]\label{th2}
For any  positive integer $N$,
there exists a simple 1-planar graph with minimum degree 6 and $n\geq N$ vertices in which
each matching is of size at most
$\frac{3}{7}n+\frac{4}{7}$.
\end{theo}

The authors in \cite{BW} suspected  that  this bound in Theorem 2 is tight and posed the following conjecture.

\begin{conj}[\cite{BW}]\label{con1}
Any  1-planar graph with minimum degree 6 and $n\geq N$ vertices has a matching of size at least $\frac{3}{7}n+O(1)$.
\end{conj}

As for  the minimum degree $\delta=7$, both papers \cite{B} and \cite{BW} constructively gave the following result.

\begin{theo}[\cite{B,BW}]\label{th3}
For any  $N$, there exists a simple 1-planar graph with minimum degree 7 and $n\geq N$ vertices for which any matching has size at most $\frac{11}{23}n+\frac{12}{23}$.
\end{theo}

%\vskip 0.2cm
  Similarly, the authors in the  papers \cite{B} and \cite{BW} wondered  whether this bound in Theorem 3 is tight, but this remains as an open problem.

%\vskip 0.2cm

In this paper we confirm  Conjecture 1 above, and have the following result.

\begin{theo}\label{th4}
Any simple 1-planar graph with minimum degree 6 and $n \geq 36$  vertices has a matching of size at least $\frac{3}{7}n+\frac{4}{7}$, and this lower bound is tight.
\end{theo}

 The paper is organized as follows. In Section~\ref{sec2} we explain some terminology and notations, and in Section~\ref{sec3} we provide some lemmas. The proof of Theorem~\ref{th4} is given in  Section~\ref{sec4}. Some problems worthy of further study are presented  in Section~\ref{sec5}.

%    \vskip 0.4cm

\section{Terminology and notation
\label{sec2}}

All graphs considered here are simple, and possibly disconnected.
Let $V(G)$ and $E(G)$  denote
the vertex set and edge set of a graph $G$, \rered{respectively}.

Let $G$ be a graph.
The {\it degree} of a vertex $v$ in  $G$, denoted by $deg_{G}(v)$,  is the number of edges of $G$  incident to $v$ (a loop edge is  counted twice),
and the {\it minimum degree} of $G$ is
denoted by $\delta (G)$.
A {\it cycle} of $G$ is  a closed trail whose origin and internal vertices are distinct.
For any subset $A\subseteq V(G)$,
if $A\ne V(G)$,
let $G\backslash A$ be the graph
obtained from $G$ by deleting
all vertices in $A$ together with their incident edges,
and if $A\ne \emptyset$,
the subgraph of $G$ induced by $A$,
denoted by $G[A]$,
is the graph $G\backslash (V(G)\setminus A)$.

A component $F$ in a graph is called
a  {\it $k$-vertex-component}
if $F$ has exactly  $k$ vertices,
and moreover  $F$ is called an {\it odd component}
({\it \rered{even component}})
if $k$ is odd (even).
For  $S\subseteq V(G)$,
let  $odd(G\backslash S)$ denote the number of odd components of $G\backslash S$.

For any two disjoint  vertex subsets $A$ and $B$
of a graph $G$, let $E_{G}(A, B)$ denote
the set of edges in $G$
which have one end vertex in $A$ and the other in $B$,
and for $v\in V(G)$,
let $N_G(v, A)$ be the set of vertices in $A$
which are adjacent to $v$.
Write ``$G\cong H$" when graphs $G$ and $H$ are  isomorphic.

Let $D$ be a drawing of $G$.
An edge $e$ of $G$ is called  {\it clean} under $D$
if it does not cross with any other edge
under the drawing $D$;
a cycle $C$ of $G$ is called
{\it clean}  if each edge on  $C$ is clean under $D$.

Let $H$ be a subgraph of $G$ with a drawing $D$.  The subdrawing $D|_H$ of $H$ induced by $D$ is called a {\it restricted drawing} of $D$.

For any drawing $D$,
the {\it associated plane graph} $\cro{D}$  of $D$
is the plane graph that is obtained from $D$ by turning all crossings
of $D$ into new vertices of degree four.
A cycle $C$  of $\cro{D}$ partitions the plane into two open regions,
the bounded one (i.e., the {\it interior } of $C$)
and the  unbounded one (i.e., the {\it exterior} of $C$).
We denote  by $int_{\cro{D}}(C)$ and  $ext_{\cro{D}}(C)$  the interior and exterior of $C$, respectively.

For other terminology and notations not defined here we refer to \cite{JAB}.

%\vskip 0.7cm

\section{Preliminary results
%Some lemmas
\label{sec3}}

%\vskip 0.3cm

 In this section we give some lemmas.
The first one is the well-known
Tutter-Berge Formula,
which gives a formula on the size of
a maximum matching  of a graph.

\begin{lemma}[Tutte-Berge \cite{BS}]\label{matching}
The  size of a maximum matching $M$ of a graph $G$ with $n$ vertices equals the minimum,
over all  $S\subseteq V(G)$,
of $\frac{1}{2}(n-(odd(G\backslash S)-|S|))$.
\end{lemma}

\begin{lemma}\label{1-vertex-co}
Let $G$ be a 1-planar graph with minimum degree $\delta (G)\geq 6$, and $S$ be a subset of $V(G)$ with $|S|\geq 2$.  Denote by $a_1$
the number of 1-vertex-components of $G\backslash S$.
Then $a_1+3\leq |S|$.
\end{lemma}

\proof
Let $T$ be the set of isolated vertices in
$G\backslash S$,  namely, 1-vertex-components of $G\backslash S$,
and let $H$ be the subgraph of $G$ with vertex set
$S\cup T$ and edge set $E_G(S,T)$.
As $H$ is bipartite and $1$-planar,
by a result on the maximum size
of a bipartite $1$-planar graph
due to Karpov~\cite{Kar}, we have
$$
|E(H)|\le 3|V(H)|-8=3(|S|+|T|)-8=3|S|+3a_1-8.
$$
Since $\delta(G)\ge 6$,
%$|E(H)|=\sum\limits_{u\in T}d_G(u)\ge 6|T|=6a_1$.Thus,
$$
6a_1=6|T|\le \sum_{u\in T}deg_G(u)=|E(H)|\le 3|S|+3a_1-8,
$$
implying that
$a_1+\frac 83\le |S|$.
As both $a_1$ and $|S|$ are integers, the result holds.
\proofend

\vskip 0.2cm

The following result in~\cite{YZD} gives an upper bound
of the size of a bipartite 1-planar graph
$G(X,Y;E)$ in terms of $|X|$ and $|Y|$.

\vskip 0.2cm

\begin{lemma}[\cite{YZD}]\label{bipartite}
Let $G$ be  a  bipartite  1-planar  graph (possibly disconnected) that has   partite sets of sizes $x$ and $y$   with
$2\leq x \leq y$. Then we have that  $|E(G)|\leq 2|V(G)|+4x-12$.
\end{lemma}

We say  that two 1-planar drawings of a graph are
{\it isomorphic}  if there is a
homeomorphism of the sphere  that maps one drawing to the other.

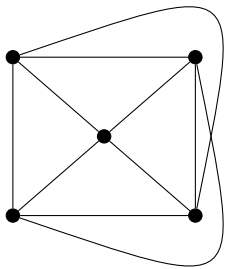
\begin{figure}[!ht]
\begin{center}
	\begin{tikzpicture}[scale=0.6]
		\begin{pgfonlayer}{nodelayer}
			\node [style=new style 0, inner sep=0pt, minimum size=5pt] (0) at (-2, 5) {};
			\node [style=new style 0, inner sep=0pt, minimum size=5pt] (1) at (-2, 1.5) {};
			\node [style=new style 0, inner sep=0pt, minimum size=5pt] (2) at (2, 5) {};
			\node [style=new style 0, inner sep=0pt, minimum size=5pt] (3) at (2, 1.5) {};
			\node [style=new style 0, inner sep=0pt, minimum size=5pt] (4) at (0, 3.25) {};
		\end{pgfonlayer}
		\begin{pgfonlayer}{edgelayer}
			\draw (0) to (2);
			\draw (2) to (3);
			\draw (0) to (1);
			\draw (1) to (3);
			\draw (0) to (4);
			\draw (4) to (2);
			\draw (4) to (1);
			\draw (4) to (3);
			\draw [bend left=60, looseness=2.65] (0) to (3);
			\draw [bend left=60, looseness=2.65] (2) to (1);
		\end{pgfonlayer}
	\end{tikzpicture}
\caption{The unique 1-planar drawing of the complete graph  $K_5$.
\label{f1}}
\end{center}
\end{figure}

%\vskip 0.2cm
\begin{lemma}\label{isomorphism}
The complete graph $K_5$ has exactly one (up to isomorphism) 1-planar drawing as shown in Figure~\ref{f1}.
\end{lemma}

\proof
See the proof of Lemma 7 in \cite{K}.
\proofend

\vskip 0.2cm

\begin{lemma}\label{5-v-com}
Let $G$ be a  1-planar graph with minimum degree $\delta(G)=6$, $S\subseteq V(G)$, and  $F$ be a 5-vertex-component of $G\backslash S$. \rered{Then the following two statements hold:}
\begin{enumerate}
\item[(a)] $|E_{G}(V(F), S)|\geq 10$; and

\item[(b)]  if  $|E_{G}(V(F), S)|=10$ or 11, then $F\cong K_5$.
\end{enumerate}
\end{lemma}

\proof (a). \rered{Since} $|V(F)|=5$, $|E(F)|\le 10$.
By the given condition,
\begin{equation}\label{eq3-1}
30\le \sum_{u\in V(F)}deg_G(u)=2|E(F)|+|E_G(V(F),S)|.
\end{equation}
Thus, $|E_G(V(F),S)|\ge 30-2\times 10=10$.

(b). If $|E_G(V(F),S)|=10$ or $11$, then
(\ref{eq3-1})  implies that
$$
|E(F)|\ge \frac 12 (30-|E_G(V(F),S)|)\ge
\frac 12(30-11)=9.5,
$$
implying that $|E(F)|\ge 10$ and so $F\cong K_5$.
\proofend

Let $D$ be a 1-planar drawing of a graph $G$.
If $L=v_1c_1v_2c_2\cdots v_{\ell}c_{\ell}v_1$
is a cycle of the associated plane graph
$\cro{D}$,
which consists alternately  of some vertices of $G$
and crossing points of $D$, then we say
that $L$ is a  {\it barrier loop} of $D$
(see \cite{MS}).
We can extend this concept ``barrier loop"
to any cycle $L$ in $\cro{D}$
which contains some clean edges $e_1,e_2,\cdots, e_k$
such that after removing these edges, each of the remaining
sections in $L$ is either an isolated vertex
or is a path in the form $v_ic_iv_{i+1}c_{i+1}\cdots v_{r+t-1} c_{r+t-1}v_{r+t}$
consisting alternately of some vertices $v_i, v_{i+1}, \cdots, v_{i+t}$ of $G$
and  crossing points  $c_i, c_{i+1}, \cdots, c_{i+t-1}$ of $D$ for some $i$ and $t$, where the subindices are taken modulo $\ell$.

The following proposition is  obvious  from  the 1-planarity.

\begin{lemma}\label{barrier}
Let  $L$ be a barrier loop of a 1-planar drawing
$D$ of a graph $G$.
For any $u,v \in V(G)$,
if $u$ and $v$ locate in $int_{\cro{D}}(L)$ and $ext_{\cro{D}}(L)$ respectively,
%regions separated by $L$,
then $u$ and $v$ are not adjacent in $G$,
and every common neighbor of $u$ and $v$
must be on $L$.
\end{lemma}

Let $D$ be a $1$-planar drawing of $G$
which has the minimum number of crossings
and let $Q$ be a restricted drawing of $D$.
Let $\F(\cro{Q})$ denote the set of faces of $\cro{Q}$.
%Let $\R(Q)$ denote the collection of regions on the plane divided by $Q$.
%Note that $Q$ divides the plane into many regions.
For any vertex $u$ in $Q$, let $\F_{u}(\cro{Q})$
be the set of faces $F$ in $\F(\cro{Q})$  such that
either $u$ is on the boundary of $F$,
or $u$ is on the boundary of some face $F_0\in \F(\cro{Q})$
whose boundary shares an edge $e$ with
the boundary of $F$, where $e$ is neither clean under $D$  and nor  crossed in $Q$.
Then, we can prove the following conclusion.

\begin{lemma}\label{u-faces}
Let $D$ be a $1$-planar drawing of $G$ and
let $Q$ be a restricted drawing of $D$.
For any $u\in V(Q)$ and any $v\in N_G(u)\setminus V(Q)$,
$v$ must be within some face of $\F_u(\cro{Q})$.
\end{lemma}

\proof As $v\notin V(Q)$, $v$ is within some face  $F$ of $\F(\cro{Q})$.

Let $Q'$ be the restricted drawing of $D$
obtained from $Q$ by adding a curve $C$ representing edge $uv$.
As $D$ is a $1$-planar drawing,
$C$ is crossed at most once.
If $C$ is not crossed in $Q'$, then
$u$ must be on the boundary of $F$
and so $F\in \F_u(\cro{Q})$.

If $C$ is crossed once in $Q'$ with some edge $e'$ in $Q$,
then $e'$ is neither clean under $D$ nor crossed in $Q$.
In $Q$, $e'$ is on the common boundary of face $F$
and another face $F_0$ in $\F(\cro{Q})$.
As $e'$ is crossed once only,
$u$ must be on the boundary of $F_0$,  implying that
$F_0\in \F_u(\cro{Q})$.
Hence $F\in \F_u(\cro{Q})$ in this case.
\proofend

\begin{figure}[!ht]
\begin{center}
	\begin{tikzpicture}[scale=0.9]
		\begin{pgfonlayer}{nodelayer}
		\node [style=new style 0, inner sep=0pt, minimum size=5pt] (0) at (2, 2) {};
		\node [style=new style 0, inner sep=0pt, minimum size=5pt] (1) at (-2, 2) {};
		\node [style=new style 0, inner sep=0pt, minimum size=5pt] (2) at (0, 3) {};
		\node [style=new style 0, inner sep=0pt, minimum size=5pt] (3) at (0, 1) {};
		\node [style=none] (10) at (2.4, 2) {$u_0$};
		\node [style=none] (11) at (-2.4,2) {$u_1$};
		\node [style=none] (12) at (0, 3.3) {$u_2$};
		\node [style=none] (13) at (0, 0.7) {$u_3$};
		\node [style=none] (14) at (0.2, 2.2) {$c$};
		%\node [style=new style 0, inner sep=0pt, minimum size=5pt] (4) at (0, 3.25) {};
		\end{pgfonlayer}
		\begin{pgfonlayer}{edgelayer}
			\draw (0) to (1);
			\draw (1) to (2);
			\draw (1) to (3);
			\draw (2) to (3);
			%\draw [bend left=60, looseness=2.65] (0) to (3);
			%\draw [bend left=60, looseness=2.65] (2) to (1);
		\end{pgfonlayer}
	\end{tikzpicture}
	
\caption{Both $u_1u_2$ and $u_1u_3$ are clean edges.
\label{f0}}
\end{center}
\end{figure}
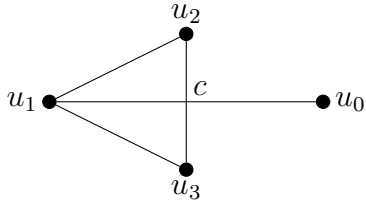

\begin{lemma}\label{clean-e}
Let $D$ be a $1$-planar drawing of $G$
which has the minimum number of crossings
and let $Q$ be a restricted drawing of $D$.
Assume that $u_1u_2u_3u_1$ is a $3$-cycle in $Q$.
If $u_2u_3$ is crossed with an edge which is incident with $u_1$,
as shown in Figure~\ref{f0},
then both $u_1u_2$ and $u_1u_3$ are clean edges  with the drawing $D$.
\end{lemma}

\proof Suppose that $u_1u_2$ is not a clean edge.
Then we redraw the edge $u_1u_2$
``most near" to one side
of the sections $u_1c$ and $cu_2$ so as to make no crossings,
contradicting to the choice of $D$ which has the  minimum  number of crossings.
\proofend

\section{The Proof of Theorem~\ref{th4}\label{sec4}
}

We only prove the former part of Theorem 4, because the tightness of the lower  bound is direct from \cite{BW}.  We first establish the following result for proving Theorem~\ref{th4}.

\begin{prop}\label{prop4-0}
Theorem~\ref{th4} is true if for every $1$-planar graph $G$
of order at least $36$
and every $S\subseteq V(G)$ with $|S|\ge 2$, the following
inequality holds:
\equ{p4-0-1}
{
3a_1+2a_3+a_5+4\leq 4|S|, \tag{*}
}
where $a_i$ is the number of $i$-vertex components
of $G\backslash S$ for $i\in \{1,3,5\}$.
\end{prop}

\proof
Let $G$ be a $1$-planar graph with $\delta(G)= 6$
and  $n\ge 36$ vertices, and
let $M(G)$ be a maximum matching of $G$.
In order to prove Theorem~\ref{th4},
i.e.,  $|M(G)|\geq \frac{3}{7}n+\frac{4}{7}$,
by Lemma \ref{matching},
it suffices to prove that,
for each subset  $S\subseteq V(G)$,
\begin{equation}\label{eq4-1}
odd(G\backslash S)-|S|\leq \frac{n-8}{7}.
\end{equation}

If $|S|=1$, because $\delta(G)=6$, each  odd component of  $G\backslash S$ has at least 7 vertices. Therefore,  $G\backslash S$ has at most  $\frac{n-1}{7}$ odd components. So,
 $$odd(G\backslash S)-|S|\leq \frac{n-1}{7}-1=\frac{n-8}{7}.$$
 Thus, (\ref{eq4-1}) holds when $|S|=1$.

 If $|S|=0$,  because  $\delta(G)=6$ and the complete graph $K_7$ is not 1-planar (see \cite{J.D}, for example),   each odd components of $G$ ($=G\backslash S$) has at least 9 vertices. Hence $G$ has at most $\frac{n}{9}$ odd components. Thus, for $n\geq 36$ we have
  $$odd(G\backslash S)-|S|\leq \frac{n}{9}\leq \frac{n-8}{7}.$$
\rered{(\ref{eq4-1}) holds as well when $|S|=0$.}

In the following we focus on the case that
$|S|\geq 2$.
%Choose arbitrarily a subset $S\subseteq V(G)$ with $|S|\geq 2$, and
%\red{\st{Assume that $D$ is a 1-planar drawing  of $G$ with minimum number of crossings.}}
For $i\geq 1$, let  $a_{2i-1}$ denote
the number of components of $G\backslash S$ with $2i-1$  vertices,
and let $a_0$ be the number of
even components of $G\backslash S$.
Noting that
$n\geq |S|+2a_0+ \sum\limits_{i\geq 1}(2i-1)a_{2i-1}$,
in order to prove (\ref{eq4-1}),
we only need  to  prove
 \begin{equation}\label{eq4-2}
 \sum\limits_{i\geq 1}a_{2i-1}-|S|\leq \frac{1}{7}\Big(|S|+2a_0+ \sum\limits_{i\geq 1}(2i-1)a_{2i-1}-8\Big).
  \end{equation}

Because $a_{2i-1}\leq \frac{2i-1}{7}a_{2i-1}$ for each $i\geq 4$,  in order to prove (\ref{eq4-2}),
we  only need to prove
\begin{equation}\label{eq4-3}
(a_1+a_3+a_5)-|S|\leq \frac{1}{7}(|S|+a_1+3a_3+5a_5-8),
\end{equation}
 namely $3a_1+2a_3+a_5+4\leq 4|S|$.

  Thus the result is proven.
\proofend

In the following, we always assume that
$G$ is a $1$-planar graph of order at least $36$
and $S$ is a subset of $V(G)$
with $|S|\ge 2$.
For each  5-vertex-component $F$ of $G\backslash S$,
$|E_{G}(V(F), S)|\geq 10$ by Lemma \ref{5-v-com} (a).
A 5-vertex component $F$ of $G\backslash S$
is called {\it bad},
if $|E_{G}(V(F), S)|=10$ or 11; otherwise, {\it good}.

The remainder of the proof  of Theorem 4  consists
of three subsections.
In Subsection~\ref{sec4-1}, we shall establish
some properties on a bad $5$-vertex component $F$ of $G\backslash S$;
in Subsection~\ref{sec4-2}, a $1$-planar bipartite graph $G^*$
will be obtained from $G$ by contracting or deleting
some edges in $G$; and in the
the last subsection, we will apply $G^*$
to show that (\ref{p4-0-1})   holds  and hence Theorem~\ref{th4} follows.

\vskip 0.4cm

\subsection{Local properties of
bad 5-vertex-components of $G\backslash S$
\label{sec4-1}}

Let $D$ be a $1$-planar drawing of $G$ such that
$D$ has the minimum number of crossings.
In this subsection,  we shall find some properties
of the local structure of a bad 5-vertex-component
of $G\backslash S$ under the drawing $D$.
Let $F$ be a bad 5-vertex-component of
$G\backslash S$ with
$V(F)=\{v_0, v_1, v_2, v_3, v_4\}$.
For each vertex $v_i\in V(F)$,
%\red{\st{denote by $N_{G}(v_i, S)$the  set of the adjacent vertices of $v$ in $S$.}}
obviously,
$\sum\limits_{i=0}^{4}|N_{G}(v_i, S)|=|E_{G}(V(F), S)|$ by the simplicity of $G$. Again, since $F$ is a bad 5-vertex-component of $G\backslash S$,
it follows from
Lemma~\ref{5-v-com} (b) that  $F\cong  K_5$.

\begin{figure}[!ht]
\begin{center}
\begin{tikzpicture}[scale=0.4]
	\begin{pgfonlayer}{nodelayer}
		\node [style=new style 0, inner sep=0pt, minimum size=5pt] (0) at (-3, 3) {};
		\node [style=new style 0, inner sep=0pt, minimum size=5pt] (1) at (-3, -3) {};
		\node [style=new style 0, inner sep=0pt, minimum size=5pt] (2) at (3, 3) {};
		\node [style=new style 0, inner sep=0pt, minimum size=5pt] (3) at (3, -3) {};
		\node [style=new style 0, inner sep=0pt, minimum size=5pt] (4) at (0, 0) {};
		\node [style=none] (10) at (-3.75, 3) {$v_1$};
		\node [style=none] (12) at (3.75, 3) {$v_2$};
		\node [style=none] (14) at (-3.75, -3) {$v_4$};
		\node [style=none] (16) at (3.75, -3) {$v_3$};
		\node [style=none] (18) at (0.85, 0) {$v_0$};
		\node [style=none] (20) at (4.3, 0) {$c$};
	\end{pgfonlayer}
	\begin{pgfonlayer}{edgelayer}
		\draw (0) to (2);
		\draw (2) to (3);
		\draw (0) to (1);
		\draw (1) to (3);
		\draw (0) to (4);
		\draw (4) to (2);
		\draw (4) to (1);
		\draw (4) to (3);
		\draw [bend left=60, looseness=2.75] (0) to (3);
		\draw [bend left=60, looseness=2.75] (2) to (1);
	\end{pgfonlayer}	
\end{tikzpicture}
\caption{The unique 1-planar drawing of  $F$ of $G\backslash S$.
\label{f2}}
\end{center}
\end{figure}
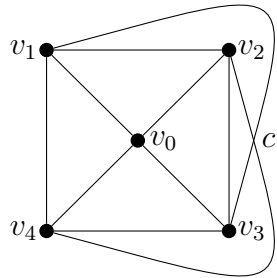

Since $D$ is a 1-planar drawing of $G$, the restricted drawing $D|_F$ is also a 1-planar drawing of $F$.
Therefore, $D|_F$ is unique up to isomorphism by
\rered{Lemma \ref{isomorphism}}.
\rered{Without} loss of generality,
in the following, \rered{we assume} that the 1-planar drawing $D|_F$ of $F$ is  depicted in Figure~\ref{f2},
where $v_1v_3$ and $v_2v_4$ are two crossed edges with the crossing point $c$.

At this time we say that  the 4-cycle $C=v_1v_2v_3v_4v_1$ is a {\it central} cycle of $F$, and $v_0$ is the {\it central} vertex  of $F$.

\rered{Next we prove} Propositions~\ref{cl1} and~\ref{cl2}.

\vskip 0.2cm

%\noindent{\bf Claim 1}. {\it
\begin{prop}\label{cl1}
For a bad 5-vertex-component $F$ of $G\backslash S$ with its central vertex $v_0$ and central cycle $C=v_1v_2v_3v_4v_1$, we have
\begin{enumerate}

\item[(a)] $|N_{G}(v_i, S)|\geq 2$ for each  $0\leq i\leq 4$;

\item[(b)] $C$  is a clean cycle under the drawing $D$;

\item[(c)]
$\rered{ \big |\bigcup\limits_{i=1}^{4}N_{G}(v_i, S)\big |}
\ge 3$;
%$|\bigcup\limits_{i=1}^{4}N_{G}(v_i, S)|\geq 3$;
and

\item[(d)] if
$\rered{ \big |\bigcup\limits_{i=1}^{4}N_{G}(v_i, S)\big |}
=3$,
%$|\bigcup\limits_{i=1}^{4}N_{G}(v_i, S)|=3$,
then
some edge $v_0v_i$, where $1\leq i\leq 4$,  is
clean  under the drawing $D$.
 \end{enumerate}
\end{prop}
%}\vskip 0.2cm

\vskip 0.2cm

\noindent
{\bf Proof}.
Let the restricted drawing $D|_F$  be as shown in Figure~\ref{f2}, where  edges $v_1v_3$ and $v_2v_4$ cross each other
at point  $c$ (see Figure~\ref{f2}).

(a) follows  immediately from the facts
that $G$ is simple with $\delta(G)=6$
and  $F\cong  K_5$.

Now prove (b).
Assume to  the contrary  that  there  exists  an edge of $C$, say $v_1v_2$,
is a crossed edge under the drawing $D$.
Then we redraw the edge $v_1v_2$
``most near" to one side
%most near along
of the sections $v_1c$ and  $cv_2$ so as to make no crossings,  contradicting to the choice of $D$ with the  minimum  number of crossings.
This proves  (b).

Then prove (c).
Obviously,  %$|\bigcup\limits_{i=1}^{4}N_{G}(v_i,S)|\geq 2$
$\rered{ \big |\bigcup\limits_{i=1}^{4}N_{G}(v_i, S)\big |}
\geq 2$
by Proposition~\ref{cl1} (a).
Assume to the contrary that
%$ |\bigcup\limits_{i=1}^{4}N_{G}(v_i, S)|\not\geq 3$.
$\rered{ \big |\bigcup\limits_{i=1}^{4}N_{G}(v_i, S)\big |}
\not\geq 3$.
Then $\bigcup\limits_{i=1}^{4}N_{G}(v_i, S)=\{x_1, x_2\}$
for two vertices $x_1,x_2$ in $S$.
Observe that
$G[\{v_0, v_1, v_2, v_3, v_4, x_1, x_2\}]$ contains a subgraph
that is isomorphic  to $K_7\backslash K_3$
(the complete graph $K_7$ by deleting three  edges of a 3-cycle).
However, it is proved in \cite{K} that  $K_7\backslash K_3$
is not 1-planar, a contradiction.  This  proves (c).

Finally prove  (d). By the assumption,
set  $\bigcup\limits_{i=1}^{4}N_{G}(v_i, S)
=\{x_1, x_2, x_3\}\subseteq S$.
For $1\le j\le 3$, let
$\ell(x_j)=|N_G(x_j,\{v_{1}, v_2, v_3, v_4\})|$.
It follows from Proposition~\ref{cl1} (a) that
$$
\sum\limits_{j=1}^{3}\ell(x_j)=\sum\limits_{i=1}^{4}|N_{G}(v_i, S)|\geq 4\times 2=8.
$$
So  there exists a vertex in $\{x_1,x_2,x_3\}$, say $x_1$,
such that  $\ell(x_1)\ge 3$,
implying that $N_G(x_1)\cap \{v_i,v_{i+1}\}\ne \emptyset$ for $i=1,2,3,4$, where $v_5$ represents the vertex $v_1$.

Note that $L_1:v_3v_4cv_3$, $L_2:v_1v_4cv_1$, $L_3:v_1v_2cv_1$ and $L_4:v_2v_3cv_2$  are barrier loops. As $v_2,v_3$ locate in $int_{\cro{D}}(L_2)$, by Lemma~\ref{barrier}, $x_1$ must be in $int_{\cro{D}}(L_2)$.
For $i=1,3,4$, as $v_i,v_{i+1}$ locate in $ext_{\cro{D}}(L_i)$, by Lemma~\ref{barrier}, $x_1$ must be in $ext_{\cro{D}}(L_i)$, where $v_5$ represents the vertex $v_1$. Thus,
$x_1$ must lie in one of the four regions  bounded by the 3-cycles  $v_1v_2v_0v_1$, $v_2v_3v_0v_2$, $v_3v_4v_0v_3$, and $v_4v_1v_0v_4$
(see Figure~\ref{f2}).
Without loss of generality,  let $x_1$ lie in the region  bounded by the 3-cycle $v_4v_1v_0v_4$.
As $\ell(x_1)\ge 3$, $x_1$ is adjacent to either $v_2$ or $v_3$,
say $v_2$, implying that $x_1v_2$ crosses $v_0v_1$.
By Lemma~\ref{clean-e}, edge $v_0v_2$ is a clean edge under the drawing $D$.

This proves (d).
\proofend

By Proposition~\ref{cl1} (c),
%$|\bigcup\limits_{i=1}^{4}N_{G}(v_i, S)|
$\rered{ \big |\bigcup\limits_{i=1}^{4}N_{G}(v_i, S)\big |}
\geq 3$.
Now we continue to describe the local
structure of $F$
under the assumption that
%$|\bigcup\limits_{i=1}^{4}N_{G}(v_i, S)|=3$
$\rered{ \big |\bigcup\limits_{i=1}^{4}N_{G}(v_i, S)\big |}
 =3$
and  $N_{G}(v_0, S)\subseteq \bigcup\limits_{i=1}^{4}N_{G}(v_i, S)$.
Now we are going to prove the following conclusion.

\begin{prop}\label{cl2}
For a bad 5-vertex component $F$ of $G\backslash S$ with
its central vertex $v_0$ and
central cycle $C: v_1v_2v_3v_4v_1$,
if $N_{G}(v_0, S)\subseteq \{x_1, x_2, x_3\}
=\bigcup\limits_{i=1}^{4}N_{G}(v_i, S)$,
then there exist two non-adjacent edges $e=v_{j_1}v_{j_2}$ and
$e'=v_{j_3}v_{j_4}$ on $C$ such that
$$
|N_G(v_{j_1})\cap N_G(v_{j_2})\cap \{x_1,x_2,x_3\}|\le 1
\quad \mbox{and}\quad
|N_G(v_{j_3})\cap N_G(v_{j_4})\cap \{x_1,x_2,x_3\}|\le 1.
$$
\end{prop}

\proof
By assumption,
$
N_{G}(v_0, S)\subseteq
\bigcup\limits_{i=1}^{4}N_{G}(v_i, S)
=\{x_1, x_2, x_3\}\subseteq S.
$
For $j=1,2,3$, let
$\ell(x_j)=|N_G(x_j, V(F))|$.
By the  simplicity of $G$,
we have that $\ell(x_j)\leq 5$
for $1\le j\le 3$,  and
$$
\sum\limits_{i=0}^{4}|N_{G}(v_i, S)|=|E_{G}(V(F), S)|=\sum\limits_{j=1}^{3}\ell(x_j).
$$
\rered{It follows  from}  Proposition~\ref{cl1} (a) that
\equ{clm2-1}
{
10=5\times 2\leq \sum\limits_{i=0}^{4}|N_{G}(v_i, S)|=\sum\limits_{j=1}^{3}\ell(x_j)=|E_{G}(V(F), S)|
\leq 11,
}
where the last inequality follows
from the assumption that $F$ is a bad $5$-vertex
component of $G\backslash S$.

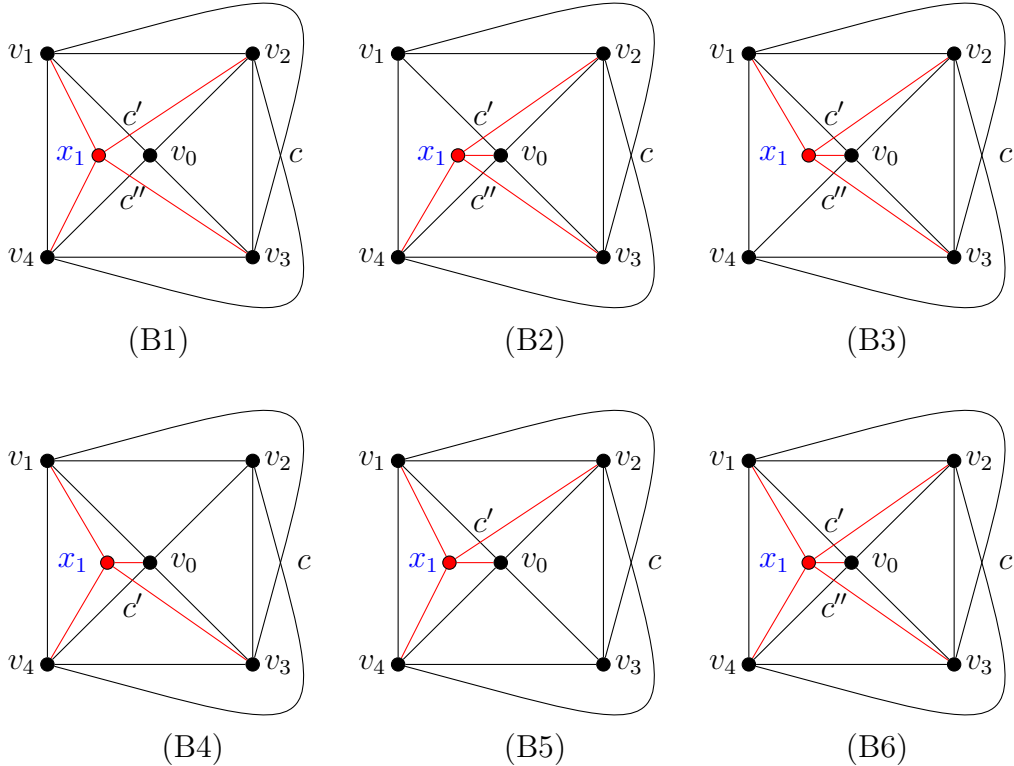
\begin{figure}[!ht]
\begin{center}
\begin{tikzpicture}[scale=0.45]
	\begin{pgfonlayer}{nodelayer}
		\node [style=new style 0, inner sep=0pt, minimum size=5pt] (0) at (1.75, -3) {};
		\node [style=new style 0, inner sep=0pt, minimum size=5pt] (1) at (1.75, 3) {};
		\node [style=new style 0, inner sep=0pt, minimum size=5pt] (2) at (7.75, -3) {};
		\node [style=new style 0, inner sep=0pt, minimum size=5pt] (3) at (7.75, 3) {};
		\node [style=new style 0, inner sep=0pt, minimum size=5pt] (4) at (4.75, 0) {};
		\node [style=none] (6) at (1, -3) {$v_4$};
		\node [style=none] (7) at (8.5, -3) {};
		\node [style=none] (8) at (8.5, -3) {$v_3$};
		\node [style=none] (10) at (1, 3) {$v_1$};
		\node [style=none] (11) at (8.5, 3) {};
		\node [style=none] (12) at (8.5, 3) {$v_2$};
		\node [style=none] (13) at (5.75, 0) {};
		\node [style=none] (14) at (5.75, 0) {$v_0$};
		\node [style=red style 1, inner sep=0pt, minimum size=5pt] (15) at (3.25, 0) {};
		\node [style=new style 0, inner sep=0pt, minimum size=5pt] (16) at (12, -3) {};
		\node [style=new style 0, inner sep=0pt, minimum size=5pt] (17) at (12, 3) {};
		\node [style=new style 0, inner sep=0pt, minimum size=5pt] (18) at (18, -3) {};
		\node [style=new style 0, inner sep=0pt, minimum size=5pt] (19) at (18, 3) {};
		\node [style=new style 0, inner sep=0pt, minimum size=5pt] (20) at (15, 0) {};
		\node [style=none] (22) at (11.25, -3) {$v_4$};
		\node [style=none] (23) at (18.75, -3) {};
		\node [style=none] (24) at (18.75, -3) {$v_3$};
		\node [style=none] (26) at (11.25, 3) {$v_1$};
		\node [style=none] (27) at (18.75, 3) {};
		\node [style=none] (28) at (18.75, 3) {$v_2$};
		\node [style=none] (29) at (16, 0) {};
		\node [style=none] (30) at (16, 0) {$v_0$};
		\node [style=new style 0, inner sep=0pt, minimum size=5pt] (32) at (22.25, -3) {};
		\node [style=new style 0, inner sep=0pt, minimum size=5pt] (33) at (22.25, 3) {};
		\node [style=new style 0, inner sep=0pt, minimum size=5pt] (34) at (28.25, -3) {};
		\node [style=new style 0, inner sep=0pt, minimum size=5pt] (35) at (28.25, 3) {};
		\node [style=new style 0, inner sep=0pt, minimum size=5pt] (36) at (25.25, 0) {};
		\node [style=none] (37) at (21.5, -3) {};
		\node [style=none] (38) at (21.5, -3) {$v_4$};
		\node [style=none] (39) at (29, -3) {};
		\node [style=none] (40) at (29, -3) {$v_3$};
		\node [style=none] (41) at (21.5, 3) {};
		\node [style=none] (42) at (21.5, 3) {$v_1$};
		\node [style=none] (43) at (29, 3) {};
		\node [style=none] (44) at (29, 3) {$v_2$};
		\node [style=none] (45) at (26.25, 0) {};
		\node [style=none] (46) at (26.25, 0) {$v_0$};
		\node [style=new style 0, inner sep=0pt, minimum size=5pt] (48) at (1.75, -15) {};
		\node [style=new style 0, inner sep=0pt, minimum size=5pt] (49) at (1.75, -9) {};
		\node [style=new style 0, inner sep=0pt, minimum size=5pt] (50) at (7.75, -15) {};
		\node [style=new style 0, inner sep=0pt, minimum size=5pt] (51) at (7.75, -9) {};
		\node [style=new style 0, inner sep=0pt, minimum size=5pt] (52) at (4.75, -12) {};
		\node [style=none] (54) at (1, -15) {$v_4$};
		\node [style=none] (55) at (8.5, -15) {};
		\node [style=none] (56) at (8.5, -15) {$v_3$};
		\node [style=none] (58) at (1, -9) {$v_1$};
		\node [style=none] (59) at (8.5, -9) {};
		\node [style=none] (60) at (8.5, -9) {$v_2$};
		\node [style=none] (61) at (5.75, -12) {};
		\node [style=none] (62) at (5.75, -12) {$v_0$};
		\node [style=new style 0, inner sep=0pt, minimum size=5pt] (64) at (12, -15) {};
		\node [style=new style 0, inner sep=0pt, minimum size=5pt] (65) at (12, -9) {};
		\node [style=new style 0, inner sep=0pt, minimum size=5pt] (66) at (18, -15) {};
		\node [style=new style 0, inner sep=0pt, minimum size=5pt] (67) at (18, -9) {};
		\node [style=new style 0, inner sep=0pt, minimum size=5pt] (68) at (15, -12) {};
		\node [style=none] (70) at (11.25, -15) {$v_4$};
		\node [style=none] (71) at (18.75, -15) {};
		\node [style=none] (72) at (18.75, -15) {$v_3$};
		\node [style=none] (74) at (11.25, -9) {$v_1$};
		\node [style=none] (75) at (18.75, -9) {};
		\node [style=none] (76) at (18.75, -9) {$v_2$};
		\node [style=none] (77) at (16, -12) {};
		\node [style=none] (78) at (16, -12) {$v_0$};
		\node [style=none] (80) at (5, -5.5) {(B1)};
		\node [style=none] (81) at (16, -5.5) {(B2)};
		\node [style=none] (82) at (26, -5.5) {(B3)};
		\node [style=none] (83) at (6, -17.5) {(B4)};
		\node [style=none] (84) at (16, -17.5) {(B5)};
		\node [style=none] (90) at (2.45, 0) {{\color{blue}{$x_1$}}};  %%original color: green
		\node [style=none] (91) at (13, 0) {{\color{blue}{$x_1$}}};%%original color: green
		\node [style=none] (92) at (23, 0) {{\color{blue}{$x_1$}}};%%original color: green
		\node [style=none] (93) at (12.75, -12) {{\color{blue}{$x_1$}}};%%original color: green
		\node [style=none] (94) at (2.5, -12) {{\color{blue}{$x_1$}}};%%original color: green
		\node [style=red style 1, inner sep=0pt, minimum size=5pt] (95) at (13.75, 0) {};
		\node [style=red style 1, inner sep=0pt, minimum size=5pt] (96) at (24, 0) {};
		\node [style=red style 1, inner sep=0pt, minimum size=5pt] (97) at (3.5, -12) {};
		\node [style=red style 1, inner sep=0pt, minimum size=5pt] (98) at (13.5, -12) {};
		\node [style=none] (99) at (9, 0) {$c$};
		\node [style=none] (100) at (4.25, 1.25) {$c'$};
		\node [style=none] (101) at (4.25, -1.25) {$c''$};
		\node [style=none] (102) at (19.25, 0) {$c$};
		\node [style=none] (103) at (14.5, 1.25) {$c'$};
		\node [style=none] (104) at (14.5, -1.25) {$c''$};
		\node [style=none] (106) at (29.75, 0) {$c$};
		\node [style=none] (107) at (24.75, 1.25) {$c'$};
		\node [style=none] (108) at (24.75, -1.25) {$c''$};
		\node [style=none] (109) at (9.25, -12) {$c$};
		\node [style=none] (110) at (4.25, -13.25) {$c'$};
		\node [style=none] (111) at (19.5, -12) {$c$};
		\node [style=none] (112) at (14.5, -10.75) {$c'$};
		\node [style=new style 0, inner sep=0pt, minimum size=5pt] (114) at (22.25, -15) {};
		\node [style=new style 0, inner sep=0pt, minimum size=5pt] (115) at (22.25, -9) {};
		\node [style=new style 0, inner sep=0pt, minimum size=5pt] (116) at (28.25, -15) {};
		\node [style=new style 0, inner sep=0pt, minimum size=5pt] (117) at (28.25, -9) {};
		\node [style=new style 0, inner sep=0pt, minimum size=5pt] (118) at (25.25, -12) {};
		\node [style=none] (119) at (21.5, -15) {};
		\node [style=none] (120) at (21.5, -15) {$v_4$};
		\node [style=none] (121) at (29, -15) {};
		\node [style=none] (122) at (29, -15) {$v_3$};
		\node [style=none] (123) at (21.5, -9) {};
		\node [style=none] (124) at (21.5, -9) {$v_1$};
		\node [style=none] (125) at (29, -9) {};
		\node [style=none] (126) at (29, -9) {$v_2$};
		\node [style=none] (127) at (26.25, -12) {};
		\node [style=none] (128) at (26.25, -12) {$v_0$};
		\node [style=none] (129) at (26, -17.5) {(B6)};
		\node [style=none] (130) at (23, -12) {{\color{blue}{$x_1$}}}; %%original color: green
		\node [style=red style 1, inner sep=0pt, minimum size=5pt] (131) at (24, -12) {};
		\node [style=none] (132) at (29.75, -12) {$c$};
		\node [style=none] (133) at (24.75, -10.75) {$c'$};
		\node [style=none] (134) at (24.75, -13.25) {$c''$};
	\end{pgfonlayer}
	\begin{pgfonlayer}{edgelayer}
		\draw (0) to (2);
		\draw (2) to (3);
		\draw (0) to (1);
		\draw (1) to (3);
		\draw (0) to (4);
		\draw (4) to (2);
		\draw (4) to (1);
		\draw (4) to (3);
		\draw [bend right=60, looseness=2.75] (0) to (3);
		\draw [bend right=60, looseness=2.75] (2) to (1);
		\draw (16) to (18);
		\draw (18) to (19);
		\draw (16) to (17);
		\draw (17) to (19);
		\draw (16) to (20);
		\draw (20) to (18);
		\draw (20) to (17);
		\draw (20) to (19);
		\draw [bend right=60, looseness=2.75] (16) to (19);
		\draw [bend right=60, looseness=2.75] (18) to (17);
		\draw (32) to (34);
		\draw (34) to (35);
		\draw (32) to (33);
		\draw (33) to (35);
		\draw (32) to (36);
		\draw (36) to (34);
		\draw (36) to (33);
		\draw (36) to (35);
		\draw [bend right=60, looseness=2.75] (32) to (35);
		\draw [bend right=60, looseness=2.75] (34) to (33);
		\draw (48) to (50);
		\draw (50) to (51);
		\draw (48) to (49);
		\draw (49) to (51);
		\draw (48) to (52);
		\draw (52) to (50);
		\draw (52) to (49);
		\draw (52) to (51);
		\draw [bend right=60, looseness=2.75] (48) to (51);
		\draw [bend right=60, looseness=2.75] (50) to (49);
		\draw (64) to (66);
		\draw (66) to (67);
		\draw (64) to (65);
		\draw (65) to (67);
		\draw (64) to (68);
		\draw (68) to (66);
		\draw (68) to (65);
		\draw (68) to (67);
		\draw [bend right=60, looseness=2.75] (64) to (67);
		\draw [bend right=60, looseness=2.75] (66) to (65);
		\draw [style=new edge style 1] (1) to (15);
		\draw [style=new edge style 1] (15) to (3);
		\draw [style=new edge style 1] (15) to (0);
		\draw [style=new edge style 1] (15) to (2);
		\draw [style=new edge style 1] (95) to (20);
		\draw [style=new edge style 1] (95) to (16);
		\draw [style=new edge style 1] (95) to (18);
		\draw [style=new edge style 1] (96) to (36);
		\draw [style=new edge style 1] (33) to (96);
		\draw [style=new edge style 1] (96) to (35);
		\draw [style=new edge style 1] (97) to (52);
		\draw [style=new edge style 1] (97) to (48);
		\draw [style=new edge style 1] (97) to (50);
		\draw [style=new edge style 1] (65) to (98);
		\draw [style=new edge style 1] (98) to (67);
		\draw [style=new edge style 1] (98) to (68);
		\draw [style=new edge style 1] (95) to (19);  %%%original style 3
		\draw [style=new edge style 1] (96) to (34);  %%%original style 3
		\draw [style=new edge style 1] (49) to (97);  %%%original style 3
		\draw [style=new edge style 1] (98) to (64);  %%%original style 3
		\draw (114) to (116);
		\draw (116) to (117);
		\draw (114) to (115);
		\draw (115) to (117);
		\draw (114) to (118);
		\draw (118) to (116);
		\draw (118) to (115);
		\draw (118) to (117);
		\draw [bend right=60, looseness=2.75] (114) to (117);
		\draw [bend right=60, looseness=2.75] (116) to (115);
		\draw [style=new edge style 1] (131) to (118);
		\draw [style=new edge style 1] (115) to (131);
		\draw [style=new edge style 1] (131) to (117);
		\draw [style=new edge style 1] (131) to (114);
		\draw [style=new edge style 1] (131) to (116);
	\end{pgfonlayer}
\end{tikzpicture}

\caption{The possible subdrawings involving $x_1$ 
%in the proof of Claim 2.
\label{f4}}
\end{center}
\end{figure}

Assume that $\ell(x_1)\ge \ell(x_i)$ for $i=2,3$.
By (\ref{clm2-1}), $\ell(x_1)\geq 4$. Note that $L_1:v_3v_4cv_3$, $L_2:v_1v_4cv_1$, $L_3:v_1v_2cv_1$ and $L_4:v_2v_3cv_2$  are barrier loops. As $v_0,v_2,v_3$ locate in $int_{\cro{D}}(L_2)$, by Lemma~\ref{barrier}, $x_1$ must be in $int_{\cro{D}}(L_2)$.
For $i=1,3,4$, as $v_0,v_i,v_{i+1}$ locate in $ext_{\cro{D}}(L_i)$, by Lemma~\ref{barrier}, $x_1$ must be in $ext_{\cro{D}}(L_i)$, where $v_5$ represents the vertex $v_1$. Thus, $x_1$ must  lie in one of the four regions bounded
by the 3-cycles $v_1v_2v_0v_1$, $v_2v_3v_0v_2$, $v_3v_4v_0v_3$, and $v_4v_1v_0v_4$
(see Figure~\ref{f2}).
Without loss of generality,
in the following we always assume
that $x_1$ lies in the region  bounded by the 3-cycle $v_4v_1v_0v_4$ (for it is completely analogous for  other cases).
Then there are six possible subdrawings (B1)-(B6),
as shown in Figure~\ref{f4},
where the first five subdrawings (B1)-(B5)
correspond to that $\ell (x_1)=4$, and the
last subdrawing (B6)  corresponds to that $\ell(x_1)=5$.

%\incase  $\ell(x_1)=4$.

We shall prove the following claims to complete the proof.

\begin{claim}\label{claim0}
For $0\le i\le 4$ and $1\le t\le 3$,
$N_G(v_i)\cap
\rered{(\{x_1,x_2,x_3\}\setminus \{x_t\})}\ne \emptyset$,
and
if $x_t\notin N_G(v_i)$, then $\{x_1,x_2,x_3\}\setminus \{x_t\}\subseteq N_G(v_i)$.
\end{claim}

\proof For $0\le i\le 4$, by the given condition
and Proposition~\ref{cl1} (a),
$$
|N_G(v_i)\cap \{x_1,x_2,x_3\}|=|N_G(v_i)\cap S|\ge 2.
$$
Thus, Claim~\ref{claim0} follows.
\proofend

\begin{claim}\label{claim1}
Subdrawings (B1), (B2) and (B3) cannot occur.
\end{claim}

%We will show that subdrawings (B1), (B2) and (B3) cannot occur.
\proof
Since (B2) and (B3) are symmetric, we consider (B1) and (B2) only.

Observe that in both (B1) and (B2),
$C_1=x_1c'v_2cv_3c''x_1$ is a barrier loop
of $D$, and $v_0$ locates in $int_{\cro{D}}(C_1)$
while $v_1$ and $v_4$ locate in $ext_{\cro{D}}(C_1)$.

In (B1),
as $x_1v_0\notin E(G)$,
by Claim~\ref{claim0},
$\{x_2,x_3\}\subseteq N_G(v_0)$.
As $v_0$ is in $int_{\cro{D}}(C_1)$,
by Lemma~\ref{barrier},
both $x_2$ and $x_3$ are in $int_{\cro{D}}(C_1)$.
Since $v_1$ locates in $ext_{\cro{D}}(C_1)$,
by Lemma~\ref{barrier} again, we have
$N_G(v_1)\cap \{x_1,x_2,x_3\}\subseteq \{x_1\}$,
%\quad \mbox{and}\quadN_G(v_1)\cap \{x_1,x_2,x_3\}\subseteq \{x_1\},$$
a contradiction to Proposition~\ref{cl1} (a).

In (B2),
as $x_1v_1\notin E(G)$,
by Claim~\ref{claim0}, %Proposition~\ref{cl1} (a),
$\{x_2,x_3\}\subseteq N_G(v_1)$.
As $v_1$ is in $ext_{\cro{D}}(C_1)$,
by Lemma~\ref{barrier},
both $x_2$ and $x_3$ are in $ext_{\cro{D}}(C_1)$.
Since $v_0$ is in $int_{\cro{D}}(C_1)$,
by Lemma~\ref{barrier} again,
we have $N_G(v_0)\cap \{x_1,x_2,x_3\}=\{x_1\}$,
a contradiction to Proposition~\ref{cl1} (a).

Hence Claim~\ref{claim1} holds.
\proofend

\begin{claim}\label{claim2}
Proposition~\ref{cl2} holds
if subdrawing (B4) or (B5) occurs.
\end{claim}

\proof
We only consider (B4) because of the symmetry.
Suppose (B4) happens.

Note that in (B4), $x_1\notin N_G(v_2)$.
By Claim~\ref{claim0},
$N_G(v_2)\cap \{x_1,x_2,x_3\}=\{x_2,x_3\}$.
%$\{x_2,x_3\}\subseteq N_G(v_2)$.
By Claim~\ref{claim0} again,
$\{x_2,x_3\}\cap N_G(v_4)\ne \emptyset$.
Assume that $x_2\in N_G(v_4)$.
Then $x_2$ is adjacent to both $v_2$ and $v_4$.
We shall show that
\equ{cl2-1}
{x_2\notin N_G(v_3)\qquad \mbox{and}\qquad
x_3\notin N_G(v_1)\cup N_G(v_4).
}

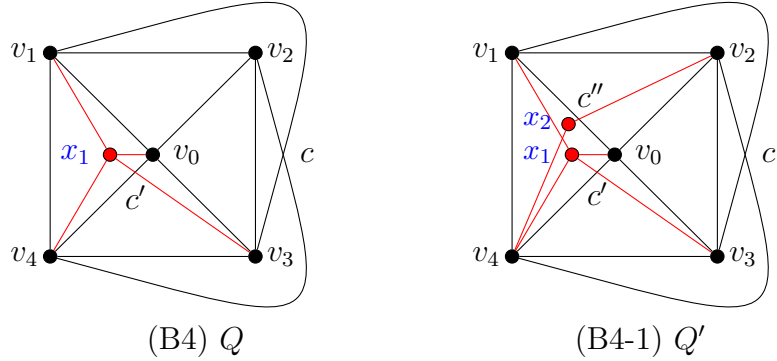
\begin{figure}[!ht]
\begin{center}
\begin{tikzpicture}[scale=0.45]
	\begin{pgfonlayer}{nodelayer}
		\node [style=new style 0, inner sep=0pt, minimum size=5pt] (48) at (1.75, -15) {};
		\node [style=new style 0, inner sep=0pt, minimum size=5pt] (49) at (1.75, -9) {};
		\node [style=new style 0, inner sep=0pt, minimum size=5pt] (50) at (7.75, -15) {};
		\node [style=new style 0, inner sep=0pt, minimum size=5pt] (51) at (7.75, -9) {};
		\node [style=new style 0, inner sep=0pt, minimum size=5pt] (52) at (4.75, -12) {};
		\node [style=none] (54) at (1, -15) {$v_4$};
		\node [style=none] (55) at (8.5, -15) {};
		\node [style=none] (56) at (8.5, -15) {$v_3$};
		\node [style=none] (58) at (1, -9) {$v_1$};
		\node [style=none] (59) at (8.5, -9) {};
		\node [style=none] (60) at (8.5, -9) {$v_2$};
		\node [style=none] (61) at (5.75, -12) {};
		\node [style=none] (62) at (5.75, -12) {$v_0$};
		\node [style=none] (83) at (6, -17.5) {(B4) $Q$};
		\node [style=none] (94) at (2.5, -12) {{\color{blue}{$x_1$}}};
		\node [style=red style 1, inner sep=0pt, minimum size=5pt] (97) at (3.5, -12) {};
		\node [style=none] (109) at (9.25, -12) {$c$};
		\node [style=none] (110) at (4.25, -13.25) {$c'$};
	\end{pgfonlayer}
		\begin{pgfonlayer}{edgelayer}
		\draw (48) to (50);
		\draw (50) to (51);
		\draw (48) to (49);
		\draw (49) to (51);
		\draw (48) to (52);
		\draw (52) to (50);
		\draw (52) to (49);
		\draw (52) to (51);
		\draw [bend right=60, looseness=2.75] (48) to (51);
		\draw [bend right=60, looseness=2.75] (50) to (49);
		\draw [style=new edge style 1] (52) to (97);
		\draw [style=new edge style 1] (97) to (48);
		\draw [style=new edge style 1] (97) to (50);
		\draw [style=new edge style 1] (49) to (97);
		%\draw [style=new edge style 1] (231) to (118);
	\end{pgfonlayer}
	%%%%%%%%%%%%%%%%%%%%%%%%%%%%B4-1 below
	
	\begin{pgfonlayer}{nodelayer}
		\node [style=new style 0, inner sep=0pt, minimum size=5pt] (114) at (15.25, -15) {};
		\node [style=new style 0, inner sep=0pt, minimum size=5pt] (115) at (15.25, -9) {};
		\node [style=new style 0, inner sep=0pt, minimum size=5pt] (116) at (21.25, -15) {};
		\node [style=new style 0, inner sep=0pt, minimum size=5pt] (117) at (21.25, -9) {};
		\node [style=new style 0, inner sep=0pt, minimum size=5pt] (118) at (18.25, -12) {};
		\node [style=none] (119) at (14.5, -15) {};
		\node [style=none] (120) at (14.5, -15) {$v_4$};
		\node [style=none] (121) at (22, -15) {};
		\node [style=none] (122) at (22, -15) {$v_3$};
		\node [style=none] (123) at (14.5, -9) {};
		\node [style=none] (124) at (14.5, -9) {$v_1$};
		\node [style=none] (125) at (22, -9) {};
		\node [style=none] (126) at (22, -9) {$v_2$};
		\node [style=none] (127) at (19.25, -12) {};
		\node [style=none] (128) at (19.25, -12) {$v_0$};
		\node [style=none] (129) at (19, -17.5) {(B4-1) $Q'$};
		\node [style=none] (130) at (16, -12) {{\color{blue}{\small $x_1$}}};
		\node [style=red style 1, inner sep=0pt, minimum size=5pt] (131) at (17, -12) {};
		\node [style=none] (230) at (16, -11) {{\color{blue}{\small $x_2$}}};
		\node [style=red style 1, inner sep=0pt, minimum size=5pt] (231) at (16.9, -11.1) {};
		%\node [style=none] (130) at (17.2, -11.5) {{\color{blue}{\small $d$}}};
		
		\node [style=none] (132) at (22.75, -12) {$c$};
		\node [style=none] (133) at (17.55, -10.25) {$c''$};
		\node [style=none] (134) at (17.75, -13.15) {$c'$};
	\end{pgfonlayer}
	\begin{pgfonlayer}{edgelayer}
		\draw (114) to (116);
		\draw (116) to (117);
		\draw (114) to (115);
		\draw (115) to (117);
		\draw (114) to (118);
		\draw (118) to (116);
		\draw (118) to (115);
		\draw (118) to (117);
		\draw [bend right=60, looseness=2.75] (114) to (117);
		\draw [bend right=60, looseness=2.75] (116) to (115);
		\draw [style=new edge style 1] (131) to (118);
		\draw [style=new edge style 1] (115) to (131);
		%\draw [style=new edge style 1] (131) to (117);
		\draw [style=new edge style 1] (131) to (114);
		\draw [style=new edge style 1] (131) to (116);
		\draw [style=new edge style 1] (231) to (114);
		%\draw [style=new edge style 1] (231) to (115);
		\draw [style=new edge style 1] (231) to (117);
		%\draw [style=new edge style 1] (231) to (118);
	\end{pgfonlayer}
\end{tikzpicture}

\caption{Restricted drawing $Q$ and the one $Q'$ obtained 
after adding a vertex $x_2$ and edges joining $x_2$ to
both $v_2$ and $v_4$ 
%the possible subdrawings  involving $x_1$ in the proof of Claim 2.
\label{f4-1}}
\end{center}
\end{figure}

As $x_2$ is adjacent to both $v_2$ and $v_4$,
by Lemma~\ref{u-faces}, $x_2$ is within a common face $F$ \rered{
of} $\F_{v_2}(\cro{Q})$ and $\F_{v_4}(\cro{Q})$,
where $Q$ is the restricted drawing of $D$ shown in
Figure~\ref{f4-1} (B4).

It can be verified that
$\F_{v_2}(\cro{Q})$ and $\F_{v_4}(\cro{Q})$
have exactly one common face, denoted by $F_1$, whose interior is $int_{\cro{D}}(C_2)$,
where $C_2$ is the cycle $v_1x_1v_0v_1$.
Thus, $x_2$ must be within $int_{\cro{D}}(C_2)$,
as shown in Figure~\ref{f4-1} (B4-1), where $c''$ is the crossing point
involving the two edges $x_2v_2$ and $v_0v_1$.
As $F_1$ does not belong to $\F_{v_3}(\cro{Q})$,
$x_2\notin N_G(v_3)$ by Lemma~\ref{u-faces}.
%Therefore, (\ref{cl2-1}) holds.

As $x_1\notin N_G(v_2)$ and $x_2\notin N_G(v_3)$,
by Claim~\ref{claim0}, $x_3\in N_G(v_2)\cap N_G(v_3)$.

Note that $C_3=v_0c''v_2cv_4c'v_0$ is a barrier loop,
and $v_3$ is in $int_{\cro{D}}(C_3)$ while $v_1$ is in $ext_{\cro{D}}(C_3)$ (See  Figure~\ref{f4-1}(B4-1)).
As $x_3\in N_G(v_3)$,
by Lemma~\ref{barrier}, we have $x_3\notin N_G(v_1)$.
Similarly, we know that
$C_4=v_0c''v_1cv_3c'v_0$ is a barrier loop,
and $v_2$ is in $int_{\cro{D}}(C_4)$
while $v_4$ is in $ext_{\cro{D}}(C_4)$ (See  Figure~\ref{f4-1}(B4-1)).
As $x_3\in N_G(v_2)$,
by Lemma~\ref{barrier}, we have $x_3\notin N_G(v_4)$.
Hence (\ref{cl2-1}) holds.

By  (\ref{cl2-1}) % (\ref{cl2-2}) and
and the fact that $x_1\notin N_G(v_2)$,
we have
$$
N_G(v_1)\cap N_G(v_2)\cap \{x_1,x_3\}
=\emptyset
\quad \mbox{and}\quad
N_G(v_3)\cap N_G(v_4)\cap \{x_2,x_3\}=\emptyset,
$$
implying that
Proposition~\ref{cl2} holds with $e=v_1v_2$ and $e'=v_3v_4$
and  Claim~\ref{claim2} is proven.
\proofend

\begin{claim}\label{claim3}
Proposition~\ref{cl2} holds
if subdrawing (B6)  occurs.
\end{claim}

\proof Assume that (B6) happens.
We first show that
for both $i=2,3$,
\equ{cl3-1}
{
\{v_1,v_3\}\not\subseteq N_G(x_i)\quad \mbox{and}\quad
\{v_2,v_4\}\not\subseteq N_G(x_i).
}

Let $C_5$ be the barrier loop $v_0c'v_2cv_4c''v_0$.
Observe that $v_1$ is in $ext_{\cro{D}}(C_5)$ and
$v_3$ is in $int_{\cro{D}}(C_5)$.
By Lemma~\ref{barrier},
$\{v_1,v_3\}\not\subseteq N_G(x_i)$ for both $i=2,3$.
Similarly, it can be proved that
$\{v_2,v_4\}\not\subseteq N_G(x_i)$ for both $i=2,3$.
Thus, (\ref{cl3-1}) holds.
It implies that each $x_i$
is adjacent to at most two consecutive vertices
on $C$ for both $i=2,3$. On the other hand,  by Claim~\ref{claim0},
\equ{cl3-3}
{
\{v_1,v_2,v_3,v_4\}\subseteq N_G(x_2)\cup N_G(x_3),
}
Hence, for $i=2,3$, $x_i$
is adjacent to exactly two consecutive vertices on $C$
and $x_2,x_3$ together are adjacent to all four vertives
on $C$. Equivalently,
%Hence, by (\ref{cl3-1}) and  (\ref{cl3-3}),
one of the two cases below happens
for some $t\in \{2,3\}$:
\begin{enumerate}
\item[(i)] $N_G(x_t)\cap \{v_1,v_2,v_3,v_4\}=\{v_1,v_2\}$
and  $N_G(x_{5-t})\cap \{v_1,v_2,v_3,v_4\}=\{v_3,v_4\}$, or
\item[(ii)] $N_G(x_t)\cap \{v_1,v_2,v_3,v_4\}=\{v_1,v_4\}$
and  $N_G(x_{5-t})\cap \{v_1,v_2,v_3,v_4\}=\{v_2,v_3\}$.
\end{enumerate}

But we can show that Case (i) above cannot happen.
Note that $C_6: x_1c'v_2cv_3c''x_1$ is a barrier loop,
and $v_0$ is in $int_D(C_6)$ while both $v_1$ and $v_4$ are in $ext_D(C_6)$.
If Case (i) happens,
then $N_G(x_t)\cap \{v_1,v_4\}\ne \emptyset$ for both $t=2,3$.
By Lemma~\ref{barrier},
$x_t\notin N_G(v_0)$ for both $t=2,3$,
a contradiction to Claim~\ref{claim0}.

Thus, Case (ii) above happens, and
$$
N_G(v_1)\cap N_G(v_2)\cap \{x_2,x_3\}=\emptyset
=N_G(v_3)\cap N_G(v_4)\cap \{x_2,x_3\},
$$
implying that Proposition~\ref{cl2} holds
with $e=v_1v_2$ and $e'=v_3v_4$.

Hence Claim~\ref{claim3} holds.
\proofend

By Claims~\ref{claim1},~\ref{claim2} and~\ref{claim3},
Proposition~\ref{cl2} is proven.
\proofend

\subsection
{To form a  bipartite 1-planar graph $G^{*}$ from $G$
\label{sec4-2}}

%\noindent{\bf 2. Making $G$ into a  bipartite 1-planar graph $G^{*}$}

\vskip 0.2cm

\noindent
In this section,  associated with the given
$1$-planar graph $G$ and $1$-planar drawing $D$ of $G$
with the minimum number of crossings,
we perform several  operations  on $G$  so as to obtain a desired bipartite 1-planar (simple) graph  $G^{*}$,
in which one bipartite set is $S$
and \rered{the other consists} of some original
vertices in $V(G)\backslash S$ and some new vertices.

Let $F$ be  a  bad 5-vertex-component of $G\backslash S$
with  $V(F)=\{v_0, v_1, v_2, v_3, v_4\}$,  its central cycle $C=v_1v_2v_3v_4v_1$ and its central vertex $v_0$.
Then we know that   $C$ is clean under the drawing $D$
by Proposition~\ref{cl1} (a).
%Let $C(F)$ denote the set of vertices in the central cycle of $F$.
Then $\rered{\big |\bigcup\limits_{i=1}^{4}N_{G}(v_i, S)\big |}\geq 3$ by Proposition~\ref{cl1}(c).

In the following, we define three operations on $F$
according to the value of
$\rered{\big |\bigcup\limits_{i=1}^{4}N_{G}(v_i, S)\big |}
%|\bigcup\limits_{i=1}^{4}N_{G}(v_i, S)|
$
and
whether $N_G(v_0, S)$ is a subset of
$\bigcup\limits_{i=1}^{4}N_{G}(v_i, S)$.

%$\bigcup_{1\le i\le 4} N_G(v_i, S)$.

%$|\bigcup\limits_{i=1}^{4}N_{G}(v_i, S)|\geq 4$ or $|\bigcup\limits_{i=1}^{4}N_{G}(v_i, S)|=3$.

{\bf Operation A }
%$\big($ %for $F$ with
(When
$\rered{\big |\bigcup\limits_{i=1}^{4}N_{G}(v_i, S)\big |}\ge 4$)
%|\bigcup\limits_{i=1}^{4}N_{G}(v_i, S)|\geq 4\big)$.
{\it Contract all the edges  of $C$  such that the four vertices of $C$  are merged into a new  vertex $v^{*}$, and delete all loops and all possible parallel edges but one for each pair of
distinct vertices which appear after edge-contraction. }

\vskip 0.2cm
Since $C$ is clean under the drawing $D$,  the edge contraction in Operation A  is performable and does not affect  the 1-planarity.   After performing  Operation A,
$F$ is transformed  into a 2-vertex graph  $F''$  with $V(F'')=\{v^{*}, v_0\}$. Moreover, we easily  know that
$v^{*}$ is still adjacent  to  each vertex  in $\bigcup\limits_{i=1}^{4}N_{G}(v_i, S)$.
% after performing  Operation A.
%deleting  multiple edges or loops that possibly appear in Operation A.
Let $\omega (F'', S)$ denote the number of edges
with one end in $V(F'')$ and  the other end in $S$
after performing Operation A.
Because $|N_{G}(v_0, S)|\geq 2$ by Proposition~\ref{cl1} (a),  we  get that
\equ{s4-2-1}
{
\omega (F'', S)=
\rered{\bigg |\bigcup\limits_{i=1}^{4}N_{G}(v_i, S)\bigg |}
+|N_{G}(v_0, S)|\geq 6.
}

If $\rered{\big |\bigcup\limits_{i=1}^{4}N_{G}(v_i, S)\big |}=3$, it then
follows  from
Proposition~\ref{cl1}(d)
that there exists a clean edge $e$ under the drawing $D$
which joins $v_0$ to some vertex in $C$.
Without loss of generality, let $e=v_0v_1$.
%such that $e$ is clean under the drawing $D$.
We shall distinguish the two cases below:
$N_{G}(v_0, S)\not\subseteq \bigcup\limits_{i=1}^{4}N_{G}(v_i, S)$
and $N_{G}(v_0, S)\subseteq \bigcup\limits_{i=1}^{4}N_{G}(v_i, S)$.

\vskip 0.2cm

{\bf Operation B} $\big($When
$|\bigcup\limits_{i=1}^{4}N_{G}(v_i, S)|=3$ and
$N_{G}(v_0, S)\not\subseteq \bigcup\limits_{i=1}^{4}N_{G}(v_i, S)\big)$.  {\it
First contract all the edges  of $C$  such that the four vertices  of $C$  are merged into a new  vertex $v^{*}$; then contract the edge $e=v_0v_1$ $(=v_0v^{*})$ such that the two end vertices of $e$ are  merged into a new  vertex $v^{**}$;  finally
delete  all loops and
all parallel edges but one for each pair of
distinct vertices which appear after edge-contraction.
%delete  all possible multiple edges or loops appearing by  edge contraction.
}
\vskip 0.2cm

Similarly, the edge contraction in Operation B   is performable and does not affect  the 1-planarity.
After performing Operation B,
$F$ is transformed  into a 1-vertex graph
$F'$ with $V(F')=\{v^{**}\}$.
Similarly, we also easily see that $v^{**}$ is still adjacent to  each vertex in $\bigcup\limits_{i=1}^{4}N_{G}(v_i, S)\cup N_{G}(v_0, S)$  after performing Operation B.
%deleting  all  multiple edges or loops that possibly  appear in .

Let $\omega(F', S)$ denote the number of edges
with one end in $V(F')$ ($=\{v^{**}\}$) and
the other end in $S$. %after performing   Operation B.
Because of the assumption  that  $\rered{\big |\bigcup\limits_{i=1}^{4}N_{G}(v_i, S)\big |}=3$ and  $N_{G}(v_0, S)\not\subseteq \bigcup\limits_{i=1}^{4}N_{G}(v_i, S)$,
we therefore  get that
\equ{s4-2-2}
{
\omega(F', S)=\rered{
\bigg|\bigcup\limits_{i=0}^{4}N_{G}(v_i, S)
%\cup N_{G}(v_0, S)
\bigg |
\ge
\bigg|\bigcup\limits_{i=1}^{4}N_{G}(v_i, S)\bigg |}+1=4.
}

If $\rered{\big |\bigcup\limits_{i=1}^{4}N_{G}(v_i, S)\big |}=3$ and $N_{G}(v_0, S)\subseteq \bigcup\limits_{i=1}^{4}N_{G}(v_i, S)$,
let $\bigcup\limits_{i=1}^{4}N_{G}(v_i, S)
=\{x_1, x_2, x_3\}$, where $x_{1}, x_2, x_3\in S$.
Then, by Proposition~\ref{cl2},
there are  two non-adjacent edges $e$ and $e'$ on $C$
such that the two end vertices of $e$ (resp. $e'$)
have at most one common neighbor $x$ (resp. $x'$)
in $\{x_1, x_2, x_3\}$.

At this time  we define  the following operation.

\vskip 0.2cm

{\bf Operation C} $\big($When  $\rered{\big |\bigcup\limits_{i=1}^{4}N_{G}(v_i, S)\big |}=3$ and $N_{G}(v_0, S)\subseteq \bigcup\limits_{i=1}^{4}N_{G}(v_i, S)\big)$.  {\it Contract the edge $e$ (\rered{resp.} $e'$)
such that its two end vertices are merged into a new vertex
$v'$ (\rered{resp.} $v''$), and then
delete all parallel edges but one for
each pair of distinct vertices which appear after edge-contraction.}

\vskip 0.2cm
Similarly, the edge contraction in Operation C can be implemented and does not change the  1-planarity.  After  finishing  Operation C,  $F$ is transformed into a 3-vertex graph  $F'''$ with $V(F''')=\{v', v^{''}, v_0\}$.
By the reasons   as  just stated   before  Operation C, we know  that performing  Operation C yields
at most two multiple edges
with one end vertex is in $\{v', v''\}$ and
the other end vertex in $S$.

Let $\omega(F^{'''}, S)$ denote  the number of edges
with one end in $V(F'')$ and the other end in $S$
after performing Operation C.
Therefore, by Proposition~\ref{cl1} (a), we get  that
\equ{s4-2-3}
{
\omega(F^{'''}, S)\ge |N_{G}(v_0, S)|+\Big(\sum\limits_{i=1}^{4}|N_{G}(v_i, S)|-2\Big)\geq 8.
}

Now we are going to construct the desired graph $G^{*}$
from $G$ by going through all the following steps.

\begin{enumerate}\setlength{\itemindent}{1.4em}
\item [{\bf Step 0}.] (deleting components).
Delete from $G\backslash S$
all even components and all $i$-vertex odd components with $i\geq 7$
(deleting a vertex must also delete its all incident edges).

\item[{\bf Step 1}.]
Perform Operation A at every bad 5-vertex-component $F$
of $G\backslash S$ if
$\rered{\big |\bigcup\limits_{i=1}^{4}N_{G}(v_i, S)\big |}$ $\geq 4$.
By the definition of
Operation A,  each $F$ is changed into a 2-vertex graph $F''$.

\item[{\bf Step 2}.]
Perform Operation B at every bad 5-vertex-component
$F$ of $G\backslash S$ if $\rered{\big|\bigcup\limits_{i=1}^{4}
N_{G}(v_i, S)\big |}$ $=3$  and $N_{G}(v_0, S)\not\subseteq \bigcup\limits_{i=1}^{4}N_{G}(v_i, S)$.
By the definition of Operation B,
$F$ is changed into a 1-vertex graph $F'$.

\item[{\bf Step 3}.] Perform Operation C at every bad 5-vertex-components $F$ of $G\backslash S$ if $\rered {\big |\bigcup\limits_{i=1}^{4}
N_{G}(v_i, S)\big |}$ $=3$ and $N_{G}(v_0, S)\subseteq \bigcup\limits_{i=1}^{4}N_{G}(v_i, S)$.
 By the definition of Operation C, $F$ is changed into a 3-vertex graph $F'''$.

\item[{\bf Step 4}.] (deleting edges)
Delete the edges in all  3-vertex-components  and
all good 5-vertex-components of $G\backslash S$;
delete the edge in every 2-vertex graph $F''$
obtained in  Step 1,
and all edges in each 3-vertex graph $F'''$
obtained by Step 3;
and delete all edges which joint two vertices in $S$.
\end{enumerate}

 \vskip 0.2cm

 After doing  all the steps above, we see that the resulting graph $G^{*}$ is a bipartite 1-planar (simple) graph with  bipartite  sets  $|S|$ and $|T|$, where $T=V(G^{*})\backslash S$.

\vskip 0.2cm

{\bf Remark}: {\it The edge contraction  in Steps 1-3  may  cause  two adjacent edges to cross with each other in the resulting 1-planar drawing of $G^{*}$. If this phenomenon  appears,  we can modify the drawing so that  this two adjacent edges are no longer crossed.}

\subsection{Applying $G^{*}$ to prove Theorem~\ref{th4}
%the inequality (*)
\label{sec4-3}}

%\vskip 0.5cm
%\noindent{\bf 3. Applying  $G^{*}$ to prove the inequality $(*)$}

\noindent
From the previous subsection,
we know that $G^{*}$ is a $1$-planar bipartite graph with
bipartite set $S$ and $T$.
Now we compute the sizes of $|T|$ and $|E(G^{*})|$.
Form the process of obtaining $G^{*}$ from $G$,
we have the following facts on $T$ and $E_{G^*}(S,T)$.

\begin{enumerate}\setlength{\itemindent}{2.2em}
\item[{\bf  Fact (1)}.]  Each 1-vertex-component $F$ of $G\backslash S$ \rered{exactly contributes}  one vertex to $T$, and
at least 6 edges to $G^{*}$. %  because $\delta (G)= 6$.

\item[{\bf  Fact (2)}.] Each 3-vertex-component $F$ of $G\backslash S$ \rered{exactly contributes} 3 vertices to $T$
and at least $12$ edges to $G^{*}$,
because
$|E_{G}(V(F), S)|\geq 3\times (\delta(G)-2)\ge 3\times 4=12$.
%Thus  $F$ contributes at least 12 edges to $G^{*}$.

\item[{\bf  Fact (3)}.] Each good 5-vertex-component $F$ of $G\backslash S$ \rered{exactly contributes} 5 vertices to $T$
and at least $12$ edges to $G^{*}$,
because $|E_{G}(V(F), S)|\geq 12$ by the definition of good 5-vertex-components of $G\backslash S$.

%Lemma \ref{5-v-com}.

\item[{\bf  Fact (4)}.]
Each bad 5-vertex-component $F$ of $G\backslash S$
with $\rered{\big |\bigcup\limits_{i=1}^{4}N_{G}(v_i, S)\big |}\geq 4$
\rered{exactly contributes} 2 vertices to $T$ and
at least 6 edges to $G^{*}$ by (\ref{s4-2-1}).

%by the statement just after Operation A we know that $F$   contributes exactly 2 vertices to $T$ and at least 6 edges to $G^{*}$.

\item[{\bf  Fact (5)}.]
Each bad 5-vertex-component $F$ of $G\backslash S$
with $\rered{ \big |\bigcup\limits_{i=1}^{4}N_{G}(v_i, S)\big |}=3$
and $N_{G}(v_0, S)\not\subseteq
\bigcup\limits_{i=1}^{4}N_{G}(v_i, S)$
%by the statement just after Operation B,
\rered{exactly contributes} one  vertex to $T$ and at least $4$ edges to $G^{*}$ by (\ref{s4-2-2}).

\item[{\bf  Fact (6)}.]
Each bad 5-vertex-component $F$ of $G\backslash S$
with $\rered{ \big |\bigcup\limits_{i=1}^{4}N_{G}(v_i, S)\big |}=3$
and $N_{G}(v_0, S)\subseteq \bigcup\limits_{i=1}^{4}N_{G}(v_i, S)$
%by the statement just after Operation C,
\rered{exactly contributes} 3 vertices to $T$ and at least 8 edges to $G^{*}$  by (\ref{s4-2-3}).
\end{enumerate}

\vskip 0.3cm

Now we are ready to prove Theorem~\ref{th4}.
%inequality (*).

\vskip 0.3cm

\noindent {\bf Proof of Theorem~\ref{th4}}:
Denote by  $a^{'}_{5}$  the number of bad 5-vertex-components $F$ of $G\backslash S$ with %$|\bigcup\limits_{i=1}^{4}N_{G}(v_i, S)|
$\rered{ \big |\bigcup\limits_{i=1}^{4}N_{G}(v_i, S)\big |}=3$ and
$N_{G}(v_0, S)\not\subseteq \bigcup\limits_{i=1}^{4}N_{G}(v_i, S)$,
$a^{''}_{5}$  the number of bad 5-vertex-components $F$ of $G\backslash S$ with %$|\bigcup\limits_{i=1}^{4}N_{G}(v_i, S)|
$\rered{ \big |\bigcup\limits_{i=1}^{4}N_{G}(v_i, S)\big |}
\geq 4$,
and
$a^{'''}_{5}$  the number of bad 5-vertex-components $F$ of $G\backslash S$ with
$\rered{ \big |\bigcup\limits_{i=1}^{4}N_{G}(v_i, S)\big |}
=3$
and $N_{G}(v_0, S)\subseteq \bigcup\limits_{i=1}^{4}N_{G}(v_i, S)$.

\vskip 0.2cm

Therefore, based on  Facts (1)-(6) above  it follows that
\begin{eqnarray}\label{s4-3-4}
|T|&=&a_1+3a_3+5(a_5-a^{'}_5-a^{''}_5-a^{'''}_5)
+a^{'}_5+2a^{''}_5+3a^{'''}_5\nonumber \\
&=&a_1+3a_3+5a_5-4a^{'}_5-3a^{''}_5-2a^{'''}_5,
\end{eqnarray}
%and
\vspace{-1.1 cm}
\begin{eqnarray}\label{s4-3-5}
|E(G^{*})|&\geq & 6a_1+12a_3+12(a_5-a^{'}_5-a^{''}_5-a^{'''}_5)+4a^{'}_5
+6a^{''}_5+8a^{'''}_5\nonumber \\
&=&6a_1+12a_3+12a_5-8a^{'}_5-6a^{''}_5-4a^{'''}_5.
\end{eqnarray}

\noindent
{\bf Case 1}.  $|S|\geq |T|+1$.

\vskip 0.2cm

Noting  $a_5\geq a^{'}_5+a^{''}_5+a^{'''}_5$,  we have
$$ |S|\geq  (a_1+3a_3+5a_5-4a^{'}_5-3a^{''}_5-2a^{'''}_5)+1\geq a_1+3a_3+a_5+1.$$
Therefore,
\equ{s4-3-9a}
{4|S|\geq 4(a_1+3a_3+a_5+1) \geq  3a_1+2a_3+a_5+4.}
Theorem~\ref{th4} follows directly from
(\ref{s4-3-9a}) and
Proposition~\ref{prop4-0}.

\vskip 0.2cm

\noindent
{\bf Case 2}. $|S|\leq |T|$.
\vskip 0.2cm

Because $G^{*}$ is a bipartite 1-planar (simple) graph with bipartite  sets  $|S|$ and $|T|$, where $|S|\geq 2$,
%no matter whether $|S|\le |T|$ or $|S|\ge |T|+1$,
from Lemma  \ref{bipartite}, we have
\begin{eqnarray}\label{s4-3-8}
|E(G^{*})|&\leq & 2(|T|+|S|)+4|S|-12\nonumber\\
&=& 2(a_1+3a_3+5a_5-4a^{'}_5-3a^{''}_5-2a^{'''}_5+|S|)+4|S|-12\nonumber\\
&=&6|S|+2a_1+6a_3+10a_5-8a^{'}_5-6a^{''}_5-4a^{'''}_5-12.
\end{eqnarray}
By (\ref{s4-3-5}) and (\ref{s4-3-8}), we have
%$|E(G^{*})|\geq  6a_1+12a_3+12a_5-8a^{'}_5-6a^{''}_5-4a^{'''}_5$, we have
$3|S|\geq 2a_1+3a_3+a_5+6$.
On the other hand, $|S|\geq a_1+3$ by Lemma \ref{1-vertex-co}.
Therefore,
\equ{s4-3-9}
{4|S|\geq  (2a_1+3a_3+a_5+6)+(a_1+3)
=3a_1+3a_3+a_5+9
> 3a_1+2a_3+a_5+4.
}
Theorem~\ref{th4} follows directly from
(\ref{s4-3-9}) and
Proposition~\ref{prop4-0}.
\proofend

\section{Further study
\label{sec5}}

For each $1$-planar graph $G$,
as $|E(G)|\le 4|V(G)|-8$, we have  $\delta(G)\le 7$.
For any positive integer $n$ and $\delta$, where $3\le \delta\le 7$,
let $\nu_{n,\delta}$ be the maximum number such that
every $1$-planar graph $G$ of order $n$ and minimum degree $\delta$
has a matching of size $\nu_{n,\delta}$.

Theorem~\ref{th1} shows that $\nu_{n,3}\ge \frac{n+12}7$ when $n\ge 7$,
$\nu_{n,4}\ge \frac{n+4}3$ when $n\ge 20$, and
$\nu_{n,5}\ge \frac{2n+3}5$ when $n\ge 21$,
while Theorem~\ref{th4} shows that $\nu_{n,6}\ge \frac{3n+4}7$ when $n\ge 36$.

Regarding $\nu_{n,7}$, the following conjecture posed by Biedl and
Wittnebel~\cite{BW} is still open.

\begin{conj}\label{con2}
There exists an integer $N$ such that
$\nu_{n,7}\ge \frac{11n+12}{23}$ when $n\ge N$.
\end{conj}

Biedl~\cite{B} studied the maximum size of matchings
in $4$-connected (resp. $5$-connected) 1-planar graphs, and obtained the following
results.

\begin{theo}[\cite{B}]
\begin{enumerate}
\item For any integer $N$,
there exists a $4$-connected $1$-planar graph $G$ with $n \ge N$ vertices
%thathas no Hamiltonian path. In particular,
in which every matching has its size at most $\frac{n+4}3$.

\item For any integer $N$,
there exists a $5$-connected $1$-planar graph $G$ with $n \ge N$ vertices
%thathas no Hamiltonian path. In particular,
in which every matching has its size at most $\frac{n-2}2$.
\end{enumerate}
\end{theo}

Biedl~\cite{B} conjectured that
every $5$-connected $1$-planar graph with $n$ vertices has a matching of
size $\frac n2 -O(1)$.
%proposed the following conjecture.
\iffalse
\begin{conj}\label{con3}
Every $5$-connected $1$-planar graph with $n$ vertices has a matching of
size $\frac n2 -O(1)$.
\end{conj}
\fi
Recently this conjecture %Conjecture~\ref{con3}
was disapproved by Huang~\cite{Huang}
who shows that for any integer $N$, there exists a $5$-connected $1$-planar graph
$G$ with $n\ge N$ vertices such that
every matching in $G$ has its size at most $\frac{n}2-\frac 38\sqrt n$.

We end this article with the following problem on the study
of maximum matchings in $t$-connected $1$-planar graphs,
where $t\le 7$.

\begin{prob}\label{prob1}
%Given positive integers $t$ and $n$, where $t\le 7$,
Let $G$ be a $t$-connected $1$-planar graph with $n$ vertices,
where $t\le 7$,
and let $M(G)$ be a maximum matching of $G$.
Study the lower bound of $|M(G)|$.
In particular,
\begin{enumerate}
\item for $t = 5$,
%find constants $a$ and $b$
%a good lower bound for $|M(G)|$ in the form
%with either $a\ge \frac 12$ or $b>\frac 38$
is there a constant $b$
such that $|M(G)| \ge \frac{n}2 - b\sqrt{n}$
when $n\ge N$ for some integer $N$?
\item  for $t = 6$ or $7$, does $G$
have a near-perfect matching, namely,
$|M(G)|=\left \lfloor \frac{n}2\right \rfloor $ when $n\ge N$
for some integer $N$?
\end{enumerate}
\end{prob}

\end{document}